
\magnification 1200
\documentstyle{amsppt}
\loadbold
\NoBlackBoxes
\def\dbend{{\manual\char127}} 
\def\d@nger{\medbreak\begingroup\clubpenalty=10000
  \def\par{\endgraf\endgroup\medbreak} \noindent\hang\hangafter=-2
  \hbox to0pt{\hskip-\hangindent\dbend\hfill}\ninepoint}
\outer\def\danger{\d@nger}
\def\dd@nger{\medbreak\begingroup\clubpenalty=10000
  \def\par{\endgraf\endgroup\medbreak} \noindent\hang\hangafter=-2
  \hbox to0pt{\hskip-\hangindent\dbend\kern1pt\dbend\hfill}\ninepoint}
\outer\def\ddanger{\dd@nger}
\def\today{\ifcase\month\or January\or February\or March\or
April\or May\or June\or July\or August\or September\or October\or
November\or December\fi \space\number\day, \number\year}

\def\qed{\hfill$\square$}

\def\ju{\vskip .8truecm plus .1truecm minus .1truecm}

\def\sju{\vskip .4truecm plus .1truecm minus .1truecm}
\def\lra{\longrightarrow}

\def\stm{\setminus}

\def\se{\subseteq}



\def\NN{\Bbb N}
\def\ZZ{\Bbb Z}

\def\al{\alpha}

\def\ka{\kappa}

\def\ov{\overline}

\def\ese#1{\buildrel{#1}\over\lra}

\def\pextcn#1#2#3#4{\hbox{\rm Pext}_{#1}^{#2}({#3},{#4})}

\def\rmod{R\bold{Mod}}

\def\noi{\noindent}

\font\tencyr=wncyr10
\def\cyr{\tencyr\cyracc}

\pageno=61
\nologo

\def\qed{\hfill$\square$}
\def\noi{\noindent}

\hsize 13.5truecm
\vsize 20truecm
\parindent=18pt
\document

\rightheadtext\nofrills{ALGEBRAIC COMPACTNESS OF $\prod \bold M_{\bold \al}/\bigoplus \bold M_{\bold\al}$}
\leftheadtext\nofrills{R.~M. Dimitri\'c}
\noi{\bf International Journal of Pure and Applied Mathematics
}
\hrule
\vskip 0.1truecm
\noi{\bf Volume 14\quad No.1\quad 2004, 61-66}

\vskip 2.5truecm
\ju
\centerline{\bf ALGEBRAIC COMPACTNESS OF $\prod \bold M_{\bold \al}/\bigoplus \bold M_{\bold\al}$}
\footnote[]{Received:\quad  April 5, 2004\qquad  \hfill \copyright 2004 Academic Publications}
\vskip 0.8truecm
\centerline{Radoslav M. Dimitri\'c}
\vskip 0.3truecm
\centerline{Texas A\&M University}
\centerline{PO Box 1675, Galveston, TX 77553, USA}
\centerline{e-mail: dimitric\@tamug.edu}
\vskip 2.3truecm

\noi{\bf Abstract:} In this note, we are working within the category $\rmod$ of
(unitary, left) $R$-modules, where $R$ is a {\bf countable} ring.
It is well known (see e.g. Kie\l pi\'nski \& Simson [5], Theorem 2.2) 
that the latter condition implies 
that  the (left) pure global dimension of $R$ is at most 1.
Given an infinite index set $A$, and
a family $M_\al\in\rmod$, $\al\in A$ we are concerned with the
conditions as to when the $R$-module
$$\prod/\coprod=\prod_{\al\in A}M_\al/\bigoplus_{\al\in A}M_\al$$
is or is not algebraically compact. There are a number of special
results regarding this question and this note is meant to be an
addition to and a generalization of the set of these results.
Whether the module in the title is algebraically compact or not depends on
the numbers of algebraically compact and non-compact modules among the components
$M_\al$.
\vskip 0.3truecm
\noi{\bf AMS Subject Classification:} 16D10, 16D80, 13C13\newline
\noi{\bf Key Words:} Algebraically compact, product mod direct sum of modules,
reduced product of modules, pure global dimension 1, countable rings
\vfill
\eject
\noi Given an (infinite) cardinal $\ka$, an $R$-module $M$ is $\ka$-{\it compact},
if, every system of $\leq\ka$ linear equations  over $M$ (with unknowns $x_j$
and almost all $r_{ij}=0$):
$$\sum_{j\in J}r_{ij}x_j=m_i\in M,\quad i\in I,\quad r_{ij}\in
R,\quad |I|, |J|\leq\ka\eqno (1)
$$ 
 has a solution in $M$ whenever all its finite subsystems
have solutions (in $M$). 
A module is ({\it algebraically}) {\it compact} if it is $\ka$-compact, for every
cardinal $\ka$.
It is well-known that if 
 $M\in\rmod$ is $\ka$-compact, for some $\ka\geq
|R|$, then $M$ is algebraically compact. 
Algebraic compactness of $M$ is
equivalent to pure injectivity and this in turn is equivalent to
$\pextcn{R}{1}{X}{M}=0$, for every $X\in\rmod$.

Recall that $\prod/\coprod$ is a special case of a more general
construction of the reduced product $\prod M_\al/\Cal F$, where
$\Cal F$ is the cofinite filter on $A$. Given a subset $B\se A$,
then $\Cal F\cap B$ and $\Cal F\cap (A\stm B)$ are cofinite
filters on $B$ and on $A\stm B$ respectively, if $\Cal F$ is the
cofinite filter on $A$. One can now easily prove the following
isomorphism (alternatively use Theorem 1.10 in [2]):
$$\prod_{\al\in A}M_\al/\bigoplus_{\al\in A}M_\al\cong
\prod_{\al\in B}M_\al/\bigoplus_{\al\in B}M_\al\times
\prod_{\al\in A\stm B}M_\al/\bigoplus_{\al\in A\stm B}M_\al. \eqno (2)$$

The proof of the following  result is straightforward, since
it uses a powerful classical result of Mycielski.
           
\proclaim{Proposition 1} {\sl For every countable index set $B$,
$$\prod/\coprod=\prod_{\al\in B}M_\al/\bigoplus_{\al\in B}M_\al$$ 
is an algebraically compact $R$-module.}
\endproclaim
{\sl Proof.} Since $B$ is countable, there is a countable family of cofinite
subsets of $B$ with empty intersection. By a classical result of Mycielski [6], Theorem 1), 
$\prod/\coprod$
is $\aleph_0$-compact. This is equivalent to its algebraic compactness, since
the rings we consider here are countable. \qed

Note that this result need not hold true, if $R$ is uncountable. For instance,
if $K$ is a field and $R=K[[X,Y]]$ is the two-variable power series algebra, then
$R^\NN/R^{(\NN)}$ is not algebraically compact (see [4], Theorem 8.42).

\proclaim{Lemma 2}
{\sl Assume that pure global dimension of $R$ is $\leq 1$. If \newline 
$\bold E:\quad 0\lra A\ese{*}B\lra C\lra 0$ is a pure exact sequence and $B$ is
pure injective, then $C$ is likewise pure injective (algebraically compact).} 
\endproclaim

{\sl Proof.} Given an arbitrary $X\in\rmod$, the segment of the $\pextcn{R}{1}{X}{\bold E}$
exact sequence we are interested in is as follows:
$\dots\lra\pextcn{R}{1}{X}{B}\lra\pextcn{R}{1}{X}{C}\lra\pextcn{R}{2}{X}{A}\lra\dots$ 
Since $puregld\,\, R\leq 1$ we have $\pextcn{R}{2}{X}{A}=0$. Since $B$ is pure
injective, we have $\pextcn{R}{1}{X}{B}=0$. These facts now force
$\pextcn{R}{1}{X}{C}=0$, i.e. $C$ is pure injective. \qed

\proclaim{Proposition 3} {\sl Let $puregld\,\, R\leq 1$ and let $A$ be an arbitrary (infinite)
index set; if every $M_\al$, $\al\in A$ is
algebraically compact, then 
$\prod_{\al\in A}M_\al/\bigoplus_{\al\in A}M_\al$ is algebraically
compact.  }
\endproclaim

{\sl Proof.} It is well known that $\coprod=\bigoplus_{\al\in A}M_\al$ is a
pure submodule of $\prod=\prod_{\al\in A}M_\al$ and that $\prod$ is
algebraically compact iff all the components $M_\al$ are algebraically compact.
Appeal to Lemma 2 completes the proof. \qed

\proclaim{Theorem 4} {\sl Given any index set $A$,
let $B\se A$ be (at most) a countable set and $\forall\al\in B $, $M_\al$
is not algebraically compact, while $\forall\al\in A\stm B$,
$M_\al$ is algebraically compact. 
Then
$$\prod/\coprod=\prod_{\al\in A}M_\al/\bigoplus_{\al\in A}M_\al$$
is algebraically compact.}
\endproclaim

{\sl Proof.} By Proposition 1, the $R$-module
$\prod_{\al\in B}M_\al/\bigoplus_{\al\in B}M_\al$
is algebraically compact. By Proposition 3, 
$\prod_{\al\in A\stm B}M_\al/\bigoplus_{\al\in A\stm B}M_\al$
is likewise algebraically compact. 
Now use isomorphism (2) to conclude that $\prod/\coprod$ is algebraically compact.
\qed
\sju

Our main concern is the converse of Theorem 4: If $\prod/\coprod$ is algebraically
compact, can we conclude that at most countably many $M_\al$'s are not
algebraically compact?

Every linear system (1) has 
a short-hand representation $\mu\cdot \bold x=\bold m$, where $\mu=(r_{ij})_{i\in I, j\in J}$
is the corresponding row-finite matrix (call it the system matrix) and  $\bold x=(x_j)_{j\in J}$, $\bold m=(m_i)_{i\in I}$
are the corresponding column vectors.
The rows of matrix $\mu$ (which are the left hand sides of equations (1)) may
be viewed as elements of the free $R$-module $\oplus_{j\in J}Rx_j$. The cardinality
of these $R$-modules is $|R|2^{|J|}$. Thus the cardinality of the set of different
matrices $\mu$ representing (left-hand-sides) of (1) is at most $(|R|2^{|J|})^{|I|}=|R|^{|I|}2^{|J||I|}$.
For purposes of algebraic compactness, it suffices to consider only $|I|=|J|=\max(|R|,\aleph_0)$,
thus the latter cardinality is at most $\max(2^{|R|}, 2^{\aleph_0})$; for countable
rings this bound is $2^{\aleph_0}$. This is
an important fact that we use in the proof of the next result.
 
\proclaim{Proposition 5} {\sl Let $|A|>\max(2^{|R|}, 2^{\aleph_0})$ and $\forall\al\in A$, $M_\al$ is
not algebraically compact. Then $\prod/\coprod$ is not algebraically compact.}
\endproclaim

{\sl Proof.}
For every $M_\al$, $\al\in A$, there is a system of equations of type (1)
$$S_\al:\quad \sum_{j\in J}r_{ij}^\al x_j^\al=m_i^\al\in M_\al,\quad i\in I,\quad r_{ij}\in
R,\quad |I|=|J|=\max(|R|,\aleph_0)\eqno (3)$$
with the corresponding row finite system matrices $\mu_\al=(r_{ij}^\al)_{i\in I, j\in J}$
and the property that every finite subsystem is solvable, without the whole
system being solvable. By the observation on the number of different system
matrices $\mu_\al$, the number of different left hand sides of systems $S_\al$
is $\max(2^{|R|}, 2^{\aleph_0})$. By the assumption on the cardinality of $A$,
we conclude that there are $|A|$ many systems $S_\al$ with identical left hand
sides. Without loss of generality we assume this is correct for all $\al\in A$,
thus we consider systems (3) where the coefficients $r_{ij}^\al=r_{ij}$ do
not vary by coordinates $\al\in A$. This coefficient uniformity enables a passage to the
induced system in $\prod M_\al/\oplus M_\al$:
$$S:\quad \sum_{j\in J}r_{ij}\ov{(x_j^\al)_{\al\in A}}=
\ov{(m_i^\al)_{\al\in A}}\quad i\in I,\eqno (4)$$
(bars denote the classes mod $\oplus_{\al\in A}M_\al$).
Every finite subsystem of $S$ is equivalent to the set of 
coordinate finite subsystems of 
$S_\al$, for all but finitely many $\al\in A$. These have
solutions, which will be the coordinates of the solutions
of the original finite subsystem of $S$. But $S$ has no
global solution, for if $x_j=\ov{(s_j^\al)_{\al\in A}}, j\in J$ were 
global solutions of $S$, then $x_j^\al=s_j^\al, j\in J$ would provide 
global solutions of $S_\al$, for almost all $\al\in A$. This contradiction
then completes
the proof that $\prod M_\al/\oplus M_\al$ is not algebraically compact.
\qed

As we have not succeeded in extending the latter result to all infinite $|A|$,
we formulate the following 
\vfill
\eject

\proclaim{Conjecture}
{\sl If $|A|$ is an uncountable index set of cardinality $\leq 2^{|R|}$ and all
$M_\al\in \rmod$, $\al\in A$, are not algebraically compact, then $\prod/\coprod$
is not algebraically compact. If this is true then, for countable rings $R$,
 $\prod/\coprod$ is algebraically
compact if and only if all but countably many $M_\al\in\rmod$, $\al\in A$ are algebraically compact.}
\endproclaim

\noi {\bf Remarks.} 
There are strong indications the conjecture is correct: Gerstner [3] proved
that $\ZZ^A/\ZZ^{(A)}$ is algebraically compact, iff $A$ is countable.
A generalization follows for reduced powers of modules over countable rings: If $M\in\rmod$ is 
not algebraically compact, then use Lemma 1.2 in [1] to conclude 
that if $M^A/M^{(A)}$ is algebraically compact then $A$ must be countable.
For Abelian groups, Rychkov [7] proved that $\prod/\coprod$
is algebraically compact if and only if $A$ is countable. In fact, if $\Cal S$
denotes a set of system matrices with the property that for every $M\in\rmod$
that is not algebraically compact, there is a $\mu\in\Cal S$ that is a system
matrix for a system proving algebraic non-compactness of $M$, let $\frak n$
denote minimal cardinality of all such systems. Close inspection of the proof
of Proposition 1, ibid. seems to reveal that the RD-purity used there is not
essential, namely that it may be replaced by purity (a condition always satisfied for
Pr\"ufer domains). In that case, if $|A|>\max(\frak n,\aleph_0)$ and all $M_\al$,
$\al\in A$ are non-compact implies that $\prod/\coprod$ is non-compact.


\input cyracc.def
\font\tencyr=wncyr10        
\def\cyr{\tencyr\cyracc}    
\font\tencyi=wncyi10        
\def\cyi{\tencyi\cyracc}    

\ju
\centerline{\bf References}
\sju

[1]   B. Franzen, 
Algebraic compactness of filter quotients, {\it Proceedings Abelian Group Theory,
Oberwolfach, 1981}, Lecture Notes in Math., Springer-Verlag {\bf 874}(1981), 228-241.

[2] T.~E. Frayne \& A.~C. Morel \& D.~S Scott,
Reduced direct products, {\it Fundamenta Mathematicae}, {\bf 51}(1962), 195--228.

[3] O. Gerstner,
Algebraische kompaktheit bei Faktorgruppen von Gruppen ganzzahliger Abbildungen,
{\it Manuscripta math.}, {\bf 11}(1974), 103--109.

[4] C.~U. Jensen \& H. Lenzing,
{\it Model theoretic Algebra with particular emphasis on fields,
rings,  modules},  Gordon and Breach Science Publishers, New York (1989).  

[5] R. Kie\l pi\'nski \& D. Simson,
On pure homological dimension, {\it Bulletin de L'Acad. Polon. Sci, S\'e. Math.}, 
{\bf 23}(1975), No.1, 1--6.

[6] Jan Mycielski,
Some compactifications of general algebras, {\it Colloquium Mathematicum}, {\bf 13}(1964),
No.1, 1--9.

[7] S.~V. Rychkov,
On factor-group of the direct product of abelian groups modulo its direct sum,
{\it Math. Notes}, {\bf 29}(1981), No.3-4, 252--257.

[Orig: {\cyr S.~V. Rychkov,  O faktor-gruppe pryamogo proizvedeniya abelevyh grupp
po ih pryamo\u\i\,\,  summe,} {\cyi Matematicheskie zametki}, {\bf 29}(1981), No.4, 491--501.]

\enddocument